# Beal's Conjecture vs. "Positive Zero", Fight

By Angela N. Moore
Yale University
angela.moore@yale.edu

January 6, 2015

## Abstract
This article seeks to encourage a mathematical dialog regarding a possible solution to Beal's Conjecture. It breaks down one of the world's most difficult math problems into layman's terms and encourages people to question some of the most fundamental rules of mathematics. More specifically; it reinforces basic algebra/critical thinking skills, makes use of properties attributed to the number one and reanalyzes the definition of a positive integer in order to provide a potential counterexample to Beal's Conjecture.

## I. THE UNDEFEATED CHAMPION (BEAL'S CONJECTURE)

$A^x + B^y = C^z$

**Where *A*, *B*, *C*, *x*, *y*, and z** are positive integers with *x*, *y*, *z* > 2, then *A*, *B*, and *C* have a common prime factor.

## II. NOW FOR THE COUNTEREXAMPLE OF DOOM

**Pre-Fight:** Beal's Conjecture is never true when $(A^x = 1) + B^y = C^z$. This is because 1 has no prime factors.

**Let the Games Begin:** There are instances when positive 0 is not the same as zero. "Signed zero is zero with an associated sign. In ordinary arithmetic, $-0 = +0 = 0$. However, in computing, some number representations allow for the existence of two zeros, often denoted by $-0$ (negative zero) and $+0$ (positive zero)."[i] Furthermore, the same website proved that signed 0 sometimes produces different results than 0. "…the concept of signed zero runs contrary to the general assumption made in most mathematical fields (and in most mathematics courses) that negative zero is the same thing as zero. Representations that allow negative zero can be a source of errors in programs, as software developers do not realize (or may forget) that, while the two zero representations behave as equal under numeric comparisons, they are different bit patterns and yield different results in some operations…it is claimed that the inclusion of signed zero in IEEE 754 makes it much **easier to achieve numerical accuracy** in some critical problems in particular when computing with complex elementary functions." The following examples are specific instances of signed zero producing different results from zero in IEEE 754 the technical standard for floating-point computation.

And I quote:

The IEEE 754 floating point standard specifies the behavior of positive zero and negative zero under various operations…

**Arithmetic:** Multiplication and division follow their usual rules for combining signs:

- $\dfrac{-0}{|x|} = -0$ (for $x$ different from 0)
- $(-0) \cdot (-0) = +0$

Addition and subtraction are handled specially if the values could cancel:

- $x + (\pm 0) = x$ (for $x$ different from 0)
- $(-0) + (-0) = (-0) - (+0) = -0$
- $(+0) + (+0) = (+0) - (-0) = +0$
- $x - x = x + (-x) = +0$ (for any finite $x$, $-0$ when rounding toward negative)

Because of negative zero (and only because of it), the statements $z = -(x - y)$ and $z = (-x) - (-y)$, for floating-point variables $x$, $y$, and $z$, cannot be optimized to $z = y - x$.

Some other special rules:

- $\sqrt{-0} = -0$
- $\dfrac{-0}{-\infty} = +0$ (follows the sign rule for division)
- $\dfrac{|x|}{-0} = -\infty$ (for non-zero $x$, follows the sign rule for division)
- $\pm 0 \times \pm \infty = \text{NaN}$ (Not a Number or interrupt for indeterminate form)
- $\dfrac{\pm 0}{\pm 0} = \text{NaN}$

End of quotation.[ii]

Beal's Conjecture is grounded in Number Theory – a branch of mathematics that deals with the properties and relationships of integers; especially positive ones. Due to the aims and scopes of Number Theory; it is imperative to derive an accurate definition of a positive integer. This cannot be accomplished unless the properties of signed zero that were demonstrated above are taken into consideration. Furthermore, since floating point computation is used in Computational Number Theory - a branch of Number Theory that involves computing; it is even more apparent that the behaviors of signed zero in IEEE 754 warrant further examination.

**Final Match:** Signed zero can be used to represent different concepts "…signed zero echoes the mathematical analysis concept of approaching 0 from below as a one-sided limit, which may be denoted by $x \to 0^-$, $x \to 0-$, or $x \to \uparrow 0$. The notation "−0" may be used informally to denote a small negative number that has been rounded to zero. The concept of negative zero also has some theoretical applications in statistical mechanics and other disciplines."[i] One of the statistical/scientific uses for signed zero is the following: "In statistical mechanics, one sometimes uses negative temperatures to describe systems with population inversion, which can be considered to have a temperature greater than positive infinity, because the coefficient of energy in the population distribution function is −1/Temperature. In this context, a temperature of −0 is a (theoretical) temperature larger than any other negative temperature, corresponding to the (theoretical) maximum conceivable extent of population inversion, the opposite extreme to +0."[ii] In addition to representing a multitude of concepts, the same site referenced a number of ways in which zero can be distinguished from signed zero. "According to the IEEE 754 standard, negative zero and positive zero should compare as equal with the usual (numerical) comparison operators, like the == operators of C and Java. In those languages, special programming tricks may be needed to distinguish the two values: Type punning the number to an **integer type**, so as to look at the sign bit in the bit pattern; using the IEEE 754 copysign() function to copy the sign of the zero to some non-zero number; taking the reciprocal of the zero to obtain either $1/(+0) = +\infty$ or $1/(-0) = -\infty$ (if the division by zero exception is not trapped)…some programming languages may provide alternative comparison operators that do distinguish the two zeros. This is the case, for example, of the equals method in Java's Double class." It would not be possible to distinguish signed zero from zero by labeling it as an **integer type** if positive zero was never considered a positive **integer**. Additionally; if zero is never positive, it would not be an option to officially differentiate zero from signed zero when performing calculations. In conclusion Since 0 is an integer, it is possible for positive 0 to not be the same as 0, it is possible for positive zero to help achieve numerical accuracy and it is possible for positive zero to produce different results from zero; that means there are some rare instances when positive 0 could technically be considered a positive integer, due to the fact that it is **both positive and an integer**.

**The Finishing Blow:** If the formula $(A^x = 1) + B^y = C^z$ is used when the existence of positive 0 that can technically be considered a positive integer is allowed, then the following statement disproves Beal's Conjecture: $1^3 + (+0)^4 = 1^5$ (When reduced this is equivalent to 1+0=1). Moves performed (Values used): A=1 B=+0 C=1 x=3 y=4 z=5

## III. CAN YOU CHOOSE A WINNER?

Some people may say that Beal's Conjecture won this fight due to the fact that zero is NEVER included in "the positive integers." Furthermore, some people may say that it is highly implied that the numbers used in the final answer must be greater than zero. Others may argue that Positive Zero won the match due to the fact that it is never specifically stated in the question that a positive integer has to be greater than zero or part of the "official positive integers." The question only states that the integers must be positive and that ANY counterexample to the question posed is acceptable. Based on what you have read in this story and what you have learned in previous math classes, who do you think won the final match?  Regardless of the answer, the complexities of signed zero warrant thorough deliberation so that the gray area regarding what it truly means to be a "positive integer" can eventually be eliminated.

## IV. QUESTIONS TO CONSIDER

Are there ever exceptions to rules in mathematics? Are there certain rules that always remain true in mathematics? Should a question be judged on what is written or what is implied? When is it acceptable to alter a mathematics rule? If the definition of a positive integer changes, should Beal be allowed to alter his question? Is positive zero really different from zero? Does allowing the use of positive zero alter the original question?

## REFERENCES

[i] https://sites.google.com/site/sixdegreesofgottfriedleibniz2/degree-4-signed-zero

[ii] http://en.wikipedia.org/wiki/Signed_zero